\documentclass[preprint,12pt]{elsarticle}
\journal{J.}
\usepackage{dsfont}
\usepackage{amsfonts}
\usepackage{amsmath}
\usepackage{epsfig}
\usepackage{amssymb}
\DeclareMathAlphabet{\mathpzc}{OT1}{pzc}{m}{it}
\begin{document}
\bibliographystyle{elsarticle-harv}
\def\R{\mathbb{R}}
\def\C{\mathbb{C}}
\def\Z{\mathbb{Z}}
\def\N{\mathbb{N}}
\def\Q{\mathbb{Q}}
\def\D{\mathbb{D}}
\def\T{\mathbb{T}}
\def\hb{\hfil \break}
\def\ni{\noindent}
\def\i{\indent}
\def\a{\alpha}
\def\b{\beta}
\def\e{\epsilon}
\def\d{\delta}
\def\De{\Delta}
\def\g{\gamma}
\def\qq{\qquad}
\def\L{\Lambda}
\def\E{\cal E}
\def\G{\Gamma}
\def\F{\cal F}
\def\K{\cal K}
\def\A{\cal A}
\def\B{\cal B}
\def\M{\cal M}
\def\P{\cal P}
\def\Om{\Omega}
\def\om{\omega}
\def\s{\sigma}
\def\t{\theta}
\def\th{\theta}
\def\Th{\Theta}
\def\z{\zeta}
\def\p{\phi}
\def\P{\Phi}
\def\m{\mu}
\def\n{\nu}
\def\l{\lambda}
\def\Si{\Sigma}
\def\q{\quad}
\def\qq{\qquad}
\def\half{\frac{1}{2}}
\def\hb{\hfil \break}
\def\half{\frac{1}{2}}
\def\pa{\partial}
\def\r{\rho}
\newcommand{\nell}{\mathcal{N}}
\newcommand{\jell}{\mathcal{J}}
\newcommand{\dell}{\mathcal{D}}
\newcommand{\mell}{\mathcal{M}}
\newcommand{\sell}{\mathcal{S}}
\begin{frontmatter}
\title{Voronoi means, moving averages, and power series}
\author{N. H. Bingham\footnote{Department of Mathematics, Imperial College London, 180 Queens Gate, London, SW7 2BZ, UK; Email: N. Bingham@imperial.ac.uk.} and Bujar Gashi\footnote{Department of Mathematical Sciences, The University of Liverpool, Liverpool, L69 7ZL, UK; Email: Bujar.Gashi@liverpool.ac.uk.}}

\begin{abstract}
We introduce a {\it non-regular} generalisation of the N\"{o}rlund mean, and show its equivalence with a certain moving average. The Abelian and Tauberian theorems establish relations with convergent sequences and certain power series. A strong law of large numbers is also proved.
\end{abstract}

\begin{keyword}
Voronoi means; N\"{o}rlund means; Moving averages; Power series; Regular variation; LLN.
\end{keyword}

\end{frontmatter}

\numberwithin{equation}{section}
\newtheorem{proof}{Proof}
\newtheorem{definition}{Definition}
\newtheorem{theorem}{Theorem}
\newtheorem{lemma}{Lemma}
\newtheorem{proposition}{Proposition}
\newtheorem{remark}{Remark}
\newtheorem{corollary}{Corollary}
\newtheorem{assumption}{Asssumption}
\section{Introduction}
Let the real sequences $\{p_n,q_n,u_n\}_{n=0}^\infty$ with $u_n\neq 0$ for $n\geq 0$, be given. The real sequence $\{s_n\}_{n=0}^\infty$ has {\it Voronoi mean}\footnote{Voronoi was the first to introduce the summability method that is now known as the {\it N\"{o}rlund mean} in the Proceedings of the Eleventh Congress of Russian Naturalists and Scientists (in Russian), St. Petersburg, 1902, pp 60-61 (see~\cite{Har} page 91).} $s$, written $s_n\rightarrow s$ $(V,p_n,q_n,u_n)$, if
\begin{eqnarray}
t_n:=\frac{1}{u_n}\sum_{k=0}^np_{n-k}q_ks_k\rightarrow s\quad (n\rightarrow\infty).\label{t}
\end{eqnarray}

There are many known special cases of the Voronoi mean. The {\it generalised N\"{o}rlund mean} $(N,p_n,q_n)$ of Borwein~\cite{B} is the $(V,p_n,q_n, (p\ast q)_n)$ mean, with
\begin{eqnarray}
(p\ast q)_n:=\sum_{k=0}^np_{n-k}q_k.\nonumber
\end{eqnarray}
Other special cases are:\\

(a) the {\it Euler method} $E_p$ of order $p\in(0,1)$, which is the Voronoi mean with $p_n=(1-p)^n/n!$, $q_n=p^n/n!$, and $u_n=(p\ast q)_n$ (see~\cite{B});\\

(b) the {\it N\"{o}rlund mean} $(N,p_n)$, which is the $(V,p_n,1,(p\ast1)_n)$ mean, and for $k>0$ and $p_n=\Gamma(n+k)/\Gamma(n+1)\Gamma(k))$ becomes the {\it Ces\`{a}ro mean} $(C,k)$ (see, for example, \S 4.1 of~\cite{Har});\\

(c) the {\it weighted mean} or the {\it discontinuous Riesz mean} $(\overline{N},q_n)$, which is the  $(V,1,q_n,(1\ast q)_n)$ mean, with the further special cases of $q_n=1$  and $q_n=1/(n+1)$ giving the Ces\`{a}ro mean $(C,1)$ and the logarithmic mean $\ell$, respectively (see, for example, \S 3.8 of~\cite{Har});\\

(d) the {\it Jajte mean} - the summability method for the law of large numbers (LLN) in~\cite{jajte}, which is the $(V,1,q_n,u_n)$ mean with $\sum_{k=0}^nq_k/u_n$ not necessarily converging to $1$ as $n\rightarrow\infty$;\\

(e) the {\it Chow-Lai mean} - the summability method for the LLN in~\cite{CL}, which is the $(V,p_n,1,u_n)$ mean with $u_n\rightarrow\infty$ and $\sum_{n=0}^\infty p_n^2<\infty$.\\

The necessary and sufficient conditions for the $(V,p_n,q_n,u_n)$ mean to be regular  are (see, for example, Theorem 2 of~\cite{Har}):\\

(i) $\sum_{k=0}^n|p_{n-k}q_k|<K|u_n|$, with $K$ independent of $n$;\\

(ii) $p_{n-k}q_k/u_n\rightarrow 0$ as $n\rightarrow\infty$ for each $k\geq0$;\\

(iii) $\sum_{k=0}^np_{n-k}q_k/u_n\rightarrow 1$ as $n\rightarrow\infty$.\\

A consequence of condition (iii) is that a regular $(V,p_n,q_n,u_n)$ mean is {\it equivalent} to a regular $(N,p_n,q_n)$ mean. Thus, the introduction of a third sequence $u_n$ in (\ref{t}), which is an essential contribution of this paper, gains us nothing {\it unless the summability method is non-regular}. Moreover, the Jajte mean does not necessarily satisfy (iii), and the Chow-Lai mean never satisfies (iii) (see (d) and (e) above, respectively). For these reasons, we do {\it not} assume that the triple $(p_n,q_n,u_n)$ necessarily satisfies the regularity conditions (i)--(iii). The non-regular summability methods, apart from their intrinsic interest within summability theory, and far from being peripheral or pathological, are useful in a variety of contexts (see, for example, $\S$ \ref{rem}).\\

The {\it Voronoi convolution} of two sequences $p_n$ and $q_n$, denoted $(p\circ q)_n$, is defined as $(p\circ q)_0:=p_0q_0$, and for $n\geq 1$ as:
\begin{eqnarray}
(p\circ q)_n:=(p\ast q)_n-(p\ast q)_{n-1}.\nonumber
\end{eqnarray}

The definition of the Voronoi mean (\ref{t}) can now be rewritten as:
\begin{eqnarray}
t_n=\frac{1}{u_n}\sum_{k=0}^n (p\circ qs)_k\rightarrow s\quad (n\rightarrow\infty),\label{t2}
\end{eqnarray}
where $(p\circ qs)_n$ denotes the Voronoi convolution of $p_n$ and $q_ns_n$.\\

Let the non-zero function $u$ be such that $u(n):=u_n$. The sequence $s_n$ has {\it continuous} Voronoi mean $s$, written $s_n\rightarrow s$ $(V_x,p_n,q_n,u(x))$, if:
\begin{eqnarray}
t_x:=\frac{1}{u(x)}\sum_{0\leq k\leq x}(p\circ qs)_k\rightarrow s\quad (x\rightarrow\infty).\nonumber
\end{eqnarray}

The formulation (\ref{t2}) of the Voronoi mean motivates the introduction of the following summability method. Let $v_0:=u_0$ and
\begin{eqnarray}
v_n:=u_n-u_{n-1},\quad n\geq1.\nonumber
\end{eqnarray}
Also let
\begin{eqnarray}
N(x)&:=&\sum_{n=0}^\infty (p\circ qs)_nx^n,\label{Nx}\\
D(x)&:=&\sum_{n=0}^\infty v_n x^n.\label{Dx}
\end{eqnarray}
If the power series $D(x)$ has radius of convergence $R\in(0,\infty]$, then $s_n$ is summable to $s$ by the {\it Voronoi power series}, written $s_n\rightarrow s$ $(\mathcal{P},p_n,q_n,v_n)$ (or, if more appropriate, $(\mathcal{P},p_n,q_n,D(x))$ , if
\begin{eqnarray}
T(x):=\frac{N(x)}{D(x)}\rightarrow s\quad(x\rightarrow R-).\label{T}
\end{eqnarray}
Three known special cases are (see, for example,~\cite{Har}):\\

($\alpha$) the {\it Abel method} $A$, which is $(\mathcal{P},1,1,1/(1-x))$ with $R=1$;\\

($\beta$) the {\it Borel method} $B$, which is $(\mathcal{P},1,1/n!,e^x)$ with $R=\infty$;\\

($\lambda$) the {\it logarithmic method} $L$, which is $(\mathcal{P},1,1/(1+n),-\log(1-x))$ with $R=1$.\\

In~\cite{BiG1}, we introduced a certain {\it moving average} summability method, which is equivalent to the logarithmic mean $\ell$. Here we introduce its {\it generalisation} appropriate for the Voronoi mean. If the function $u$ is invertible, and $u(x)\sim u([x])$, where $[\cdot]$ denotes the integer part of $x$, then for $\lambda\in(1,\infty)$ we define
\begin{eqnarray}
w_\lambda(x):= u^{\leftarrow}(u(x)/\lambda),\nonumber
\end{eqnarray}
where $u^{\leftarrow}$ denotes the inverse function of $u$. In this case, the sequence $s_n$ has {\it Voronoi moving average} $s$, written $s_n\rightarrow s$ $(\mathcal{V}, p_n, q_n, u_n,\lambda)$, if
\begin{eqnarray}
c_n:=\frac{1}{u(n)}\sum_{w_\lambda(n)<k\leq n}(p\circ qs)_k\rightarrow (1-\lambda^{-1})s\quad (n\rightarrow\infty).\label{mov}
\end{eqnarray}
We write $s_n\rightarrow s$  $(\mathcal{V}_x, p_n, q_n, u(x),\lambda)$ if the limit is taken through a continuous variable. Two known special cases of this method are:\\

($\delta$) the {\it deferred Ces\`{a}ro mean} $(D,n/\lambda,n)$ of Agnew~\cite{Agn}, which is the $(\mathcal{V},1,1,n,\lambda)$ average;\\

($\mu$) the {\it logarithmic moving average} $\mathcal{L}(\lambda)$ of~\cite{BiG1}, which is the $(\mathcal{V},1,1/(1+n),\log n,\lambda)$ average.\\

The next section states our results on the properties of the introduced methods, the relations between them, and a law of large numbers. In $\S$ \ref{pro} we give the proofs, and conclude with some further remarks in the last section.

\section{Results}\label{res}
We begin with some necessary and sufficient conditions for the sequence $s_n$ to have a Voronoi mean. Recall $v_n:=u_n-u_{n-1}$.
\begin{theorem}\label{Conv} Let $u_n$ be a positive and monotonically increasing sequence such that $u_n\rightarrow\infty$ as $n\rightarrow\infty$. We have $s_n\rightarrow s$ $(V,p_n,q_n,u_n)$ if and only if
\begin{eqnarray}
(p\circ qs)_n=v_na_n+b_n,\label{conv_cond}
\end{eqnarray}
where $a_n\rightarrow s$ as $n\rightarrow\infty$ and $\sum_{n=0}^\infty b_n/u_n$ converges.
\end{theorem}
This is a generalisation of Theorem 6.5 of Bingham and Goldie~\cite{BinG2}, which was established for the Ces\`{a}ro mean $(C,1)$. In~\cite{BiG4}, we obtain an analogous result for integrals.\\

The following is a {\it limitation theorem} for the Voronoi means, and is a generalisation of Theorem 13 of Hardy~\cite{Har} for the $(\overline{N},q_n)$ mean.
\begin{theorem}\label{lim} Let $u_n/u_{n-1}=O(1)$. If $s_n\rightarrow s$ $(V,p_n,q_n,u_n)$, then
\begin{eqnarray}
(p\circ qs)_n=sv_n+o(u_n).\nonumber
\end{eqnarray}
\end{theorem}

The ordinary convergence $s_n\rightarrow s$ as $n\rightarrow\infty$, written $s_n\rightarrow s$ $(\Omega)$, {\it always} implies the summability of $s_n$ by a regular method. This is no longer the case if the summability method is non-regular. Our next result gives some necessary and sufficient conditions for $(\Omega)\Rightarrow (V,p_n,q_n,u_n)$. We also give conditions for the converse implication $(V,p_n,q_n,u_n)\Rightarrow (\Omega)$; this is a generalisation of Theorem 2.1 of M\'{o}ritz and Stadtm\"{u}ller~\cite{MS}, which was established for $(\overline{N},q_n)\Rightarrow(\Omega)$.\\

Let $m_0:=q_0$, and
\begin{eqnarray}
m_n:=q_n-q_{n-1},\quad n\geq 1.\nonumber
\end{eqnarray}
If $v_n$ is positive, non-increasing, and $u_n\rightarrow\infty$ as $n\rightarrow\infty$, then it can be shown that there exists a real sequence $\{h_n\}_{n=0}^\infty$ such that
\begin{eqnarray}
q_ns_n=\sum_{k=0}^nh_{n-k}u_kt_k,\quad n\geq 0\label{inv}
\end{eqnarray}
(see, for example, Ishiguro~\cite{Ish}). Following M\'{o}ritz and Stadtm\"{u}ller~\cite{MS}, we write
\begin{eqnarray}
U_q&:=&\left\{\alpha:\mathds{N}_0\rightarrow \mathds{N}_0\mid \lim_{n\rightarrow\infty}\alpha(n)\rightarrow\infty\quad\mbox{and}\quad\liminf_{n\rightarrow\infty}\frac{q_{\alpha(n)}}{q_n}>1\right\},\nonumber\\
\nonumber\\
L_q&:=&\left\{\beta:\mathds{N}_0\rightarrow \mathds{N}_0\mid \lim_{n\rightarrow\infty}\beta(n)\rightarrow\infty\quad\mbox{and}\quad\liminf_{n\rightarrow\infty}\frac{q_n}{q_{\beta(n)}}>1\right\},\nonumber\\
\nonumber\\
U_u&:=&\left\{\gamma:\mathds{N}_0\rightarrow \mathds{N}_0\mid \lim_{n\rightarrow\infty}\gamma(n)\rightarrow\infty\quad\mbox{and}\quad\liminf_{n\rightarrow\infty}\frac{u_{\gamma(n)}}{u_n}>1\right\},\nonumber\\
\nonumber\\
L_u&:=&\left\{\theta:\mathds{N}_0\rightarrow \mathds{N}_0\mid \lim_{n\rightarrow\infty}\theta(n)\rightarrow\infty\quad\mbox{and}\quad\liminf_{n\rightarrow\infty}\frac{u_n}{u_{\theta(n)}}>1\right\}.\nonumber\\
\nonumber
\end{eqnarray}
\begin{theorem}\label{oco} (i) Let $s_n\rightarrow s$ $(\Omega)$. Also let: $q_n\neq 0$ for $n\geq 0$; (\ref{inv}) hold for some sequence $h_n$; and $U_q$ and $L_q$ be non-empty. Then the necessary and sufficient conditions for $s_n\rightarrow s$ $(V,p_n,q_n,u_n)$ are:
\begin{eqnarray}
\sup_{\alpha\in U_q}\liminf_{n\rightarrow\infty}\frac{1}{q_{\alpha(n)}-q_n}\sum_{k=n+1}^{\alpha(n)}[(h\circ ut)_k-m_kt_n]&\geq&0,\label{con1}\\
\nonumber\\
\sup_{\beta\in L_q}\liminf_{n\rightarrow\infty}\frac{1}{q_n-q_{\beta(n)}}\sum_{k=\beta(n)+1}^{n}[m_kt_n-(h\circ ut)_k]&\geq&0.\label{con2}
\end{eqnarray}
(ii) Let $s_n\rightarrow s$ $(V,p_n,q_n,u_n)$. Also let $U_u$ and $L_u$ be non-empty. Then the necessary and sufficient conditions for $s_n\rightarrow s$ $(\Omega)$ are:
\begin{eqnarray}
\sup_{\gamma\in U_u}\liminf_{n\rightarrow\infty}\frac{1}{p_{\gamma(n)}-p_n}\sum_{k=n+1}^{\gamma(n)}[(p\circ qs)_k-v_ks_n]&\geq&0,\label{con3}\\
\nonumber\\
\sup_{\theta\in L_u}\liminf_{n\rightarrow\infty}\frac{1}{p_n-p_{\theta(n)}}\sum_{k=\theta(n)+1}^{n}[v_ks_n-(p\circ qs)_k]&\geq&0.\label{con4}
\end{eqnarray}
\end{theorem}
We refer to the {\it Tauberian conditions} (\ref{con1})-(\ref{con4}) as $(TCO)$, and they are {\it best possible} for the following equivalence.
\begin{corollary} Let: $q_n\neq 0$ for $n\geq 0$; (\ref{inv}) hold for some sequence $h_n$; and $U_q$, $L_q$, $U_u$, $L_u$, be nonempty. If $(TCO)$ holds, then $(\Omega)\Leftrightarrow(V,p_n,q_n,u_n)$.
\end{corollary}

There are many {\it inclusion} and {\it equivalence} theorems for $(N,p_n)$, $(\overline{N},q_n)$, $(N,p_n,q_n)$, and various special cases thereof (see, for example,~\cite{BK},~\cite{das1}, ~\cite{Har},~\cite{ish},~\cite{Ish},~\cite{Kie2},~\cite{KS},~\cite{K},~\cite{KR},~\cite{US},~\cite{UST},~\cite{Tanaka},~\cite{Th},~\cite{ZB}). Of course, all such results apply to the appropriately specialised Voronoi means. The well- known result {\it Kronecker's lemma} (see, for example, page 129 of~\cite{Knopp}) is an inclusion theorem for Voronoi means:\\

\ni{\bf Theorem K.} {\it Let $\{g_n\}_{n=0}^\infty$  be any sequence of monotone increasing positive numbers. If $s_n\rightarrow s$ $(V,1,q_n,1)$, then $s_n\rightarrow 0$ $(V,1,q_ng_n,g_n)$.}\\

\ni The following is another inclusion result where the summation to $s$ by one method implies summation to $0$ by another method.

\begin{theorem}\label{k} Let $\{u_n,q_n,\tilde{u}_n,\tilde{q_n}\}_{n=0}^\infty$ be positive sequences such that: $u_{n+1}/u_{n}\rightarrow 1$ as $n\rightarrow\infty$; $\tilde{u}_n\rightarrow\infty$ as $n\rightarrow\infty$; and $u_n/q_n\rightarrow 1$ $(V,1,\tilde{q}_n,\tilde{u}_n)$. If $s_n\rightarrow s$ $(V,1,q_n,u_n)$, then $s_n\rightarrow 0$ $(V,1,\tilde{q}_n,\tilde{u}_n)$.
\end{theorem}

We now generalise the results of~\cite{BiG1}, which were established for the logarithmic mean $\ell$. Let $\Lambda$ denote the set of all functions $u$ that are invertible and $u(x)\sim u([x])$. If we write $u_n\in\Lambda$, then we mean $u_n=u(n)$ and $u\in\Lambda$.

\begin{theorem}\label{ND} If $u_n\in\Lambda$, then
\begin{eqnarray}
(V,p_n,q_n,u_n)\Leftrightarrow (\mathcal{V},p_n,q_n,u_n,\lambda)\quad\mbox{for some (all) $\lambda\in(1,\infty)$}.
\end{eqnarray}
\end{theorem}
\begin{theorem}\label{Uni} Let $u_n\in\Lambda$. If (\ref{mov}) holds for all $\lambda\in(1,\infty)$, then it holds uniformly on compact $\lambda$-sets of $(1,\infty)$.
\end{theorem}
\begin{theorem}\label{NDJM} If $u\in\Lambda$, then
\begin{eqnarray}
s_n\rightarrow s\quad (V_x,p_n,q_n,u(x))\quad\Leftrightarrow\quad s_n\rightarrow s\quad(\mathcal{V}_x,p_n,q_n,u(x),\lambda)\quad\forall\lambda>1.\nonumber
\end{eqnarray}
\end{theorem}
\begin{theorem}\label{UN} If $u\in\Lambda$ and
\begin{eqnarray}
U(x):=\sum_{0\leq k\leq x}(p\circ qs)_k.\nonumber
\end{eqnarray}
then the following statements are equivalent:\\

\ni(i) $U(x)=U_1(x)-U_2(x)$, with $U_1(x)$ satisfying
\begin{eqnarray}
\lim_{x\rightarrow\infty}\frac{U_1(x)-U_1(w_\lambda(x))}{u(x)}=s(1-\lambda^{-1}),\quad \forall \lambda>1,\nonumber
\end{eqnarray}
and $U_2(x)$ non-decreasing,\\

\ni(ii)
\begin{eqnarray}
\liminf_{\alpha\downarrow1} \limsup_{x\rightarrow\infty}\sup_{\lambda\in[1,\alpha]}\frac{U(x)-U(w_\lambda(x))}{u(x)}<\infty.\nonumber
\end{eqnarray}
\end{theorem}

\begin{corollary}\label{xn} If $u\in\Lambda$, then
\begin{eqnarray}
(V,p_n,q_n,u(n))\Leftrightarrow (V_x,p_n,q_n,u(x)).\nonumber
\end{eqnarray}
\end{corollary}

The next theorem establishes relations between Voronoi means and Voronoi power series. Some statements require the notions of {\it slowly} and {\it regularly} varying functions, for which see~\cite{BinGT}.\\

\begin{theorem}\label{t1} (i) Let $v_n>0$; $u_n\rightarrow\infty$ as $n\rightarrow\infty$; and $R\in(0,\infty)$. If $s_n\rightarrow s$ $(V,p_n,q_n,u_n)$, then $s_n\rightarrow s$ $(\mathcal{P},p_n,q_n,v_n)$.\\

\ni(ii) Let $\rho\geq-1$; $v_n$ be regularly varying of index $\rho$; $u_n\rightarrow\infty$; and $R=1$. If $(p\circ sq)_n/v_n=O_L(1)$, then $s_n\rightarrow s$ $(\mathcal{P},p_n,q_n,v_n)$ implies $s_n\rightarrow s$ $(V,p_n,q_n,u_n)$.\\

\ni(iii) Let $U(x)\geq0$; $s\geq 0$; $\rho>-1$; $\hat{U}(s):=s\int_0^\infty e^{-sx}U(x)dx$ converge for $s>0$; $\ell(x)$ be a given slowly varying function; $u(x)=x^\rho\ell(x)/\Gamma(1+\rho)$; and
\begin{eqnarray}
\frac{D(x)(1-x)(-\log x)^\rho}{\ell(-1/\log x)}\rightarrow 1\quad(x\rightarrow 1-).\nonumber
\end{eqnarray}
If $s_n\rightarrow s$ $(V_x,p_n,q_n,u(x))$, then $s_n\rightarrow s$ $(\mathcal{P},p_n,q_n,D(x))$. Conversely, $s_n\rightarrow s$ $(\mathcal{P},p_n,q_n,D(x))$ implies $s_n\rightarrow s$ $(V_x,p_n,q_n,u(x))$ if and only if
\begin{eqnarray}
\lim_{\lambda\downarrow 1}\liminf_{x\rightarrow \infty}\inf_{t\in[1,\lambda]}\frac{1}{x^\rho\ell(x)}\sum_{x<k\leq tx}(p\circ qs)_k\geq 0.\nonumber
\end{eqnarray}

\ni(iv) Let $\rho\geq-1$ and $\rho\neq 0,1,...$; $v_n$ be regularly varying of index $\rho$; $u_n\rightarrow\infty$; and $R=1$. If
\begin{eqnarray}
(p\circ sq)_n/v_n-(p\circ sq)_{n-1}/v_{n-1}=O_L(v_n/u_n),\nonumber
\end{eqnarray}
then $s_n\rightarrow s$ $(\mathcal{P},p_n,q_n,v_n)$ implies $(p\circ sq)_n/v_n\rightarrow s$.\\

\ni(v) Let $v_n>0$; $v_n=O(1/n)$; $u_n\rightarrow\infty$; and $R=1$. If
\begin{eqnarray}
(p\circ sq)_n/v_n-(p\circ sq)_{n-1}/v_{n-1}=o(v_n/u_n),\nonumber
\end{eqnarray}
then $s_n\rightarrow s$ $(\mathcal{P},p_n,q_n,v_n)$ implies $(p\circ sq)_n/v_n\rightarrow s$.
\end{theorem}
Part (i) is an {\it Abelian} result, and it is a generalization of several known special cases (see, for example,~\cite{Har},~\cite{soni}). Part (ii) is Tauberian and is a generalisation of Theorem 4.1 of~\cite{jakim} established for the $J$-method (see, for example,~\cite{Har}) and $(\overline{N},p_n)$. One can extend other closely related Tauberian results, such as those in~\cite{kratz}, in a similar way. Part (iii) contains a Tauberian result of best possible character, and it is a specialization of the Hardy-Littlewood-Karamata theorem for the Laplace-Stieltjes transform. Similar results for Abel and $L$ methods of summation appear in~\cite{Bin4} and~\cite{BiG1}, respectively. Parts (iv) and (v) are a certain generalisation of Theorem 5.3 of~\cite{jakim} and Theorem 1 of~\cite{Ish7}, respectively, which were established for the $J$-method and convergence of $s_n$. Other results of this nature (see, for example,~\cite{jakim},~\cite{Ish6},~\cite{kratz}) can be extended similarly.\\

In~\cite{jajte}, Jajte introduced a law of large numbers for the $(V,1,q_n,u_n)$ summability method. We extend his result by including equivalence relations with other summability methods. Moreover, we generalise the results of~\cite{BiG1} on the LLN of Baum-Katz type, which were obtained for the logarithmic mean $\ell$.\\

In Theorem \ref{LLN} below, we encounter infinite families of summability methods which, while by no means equivalent, {\it become equivalent in the LLN context}, to the same moment condition.  This interesting phenomenon goes back to Chow~\cite{Chow} in 1973 (Euler methods; finite variance) and Lai~\cite{Lai} in 1974 (Ces\`aro means $(C,\alpha)$, $\alpha \geq 1$; finite mean), and has been developed by, e.g. the first author (\cite{Bin5},~\cite{Bin1},~\cite{Bin10}).\\

Let $\phi:[0,\infty)\rightarrow(0,\infty)$ be such that:\\

(i) $\phi(x)$ is strictly increasing,\\

(ii) $\phi(x+1)/\phi(x)\leq c$ for some constant $c>0$,\\

(iii) for some positive constants $a$ and $b$ it holds that

\begin{eqnarray}
\phi^2(s)\int_s^\infty\frac{dx}{\phi(x)^2}\leq a s+b,\quad s>0.\nonumber
\end{eqnarray}

\begin{theorem}\label{LLN} Let $X,X_1,X_2,...$, be a sequence of $i.i.d.$ random variables, and $m_k:=\mathds{E}[X_k\mathds{1}_{\{|X_k|\leq \phi(k)\}}]$.\\

\ni(a) Let the sets $\Phi_V(\phi)$ and $\widetilde{\Phi}_V(\phi)$ be defined as:
\begin{eqnarray}
\Phi_V(\phi)&:=&\{(u_n,q_n):\quad\mbox{$u_n>0$ and increasing; $q_n>0$; and $u_n/q_n=\phi(n)$}\},\nonumber\\
\nonumber\\
\widetilde{\Phi}_V(\phi)&:=&\{(u_n,q_n)\in\Phi_V(\phi):\mbox{$v_n\geq\sigma>0$ and $v_{n+1}v_{n-1}\geq v_n^2$}\}.\nonumber
\end{eqnarray}
\ni The following four statements are equivalent:\\

$(i)$ $\mathds{E}\left[\phi^{\leftarrow}(|X|)\right] < \infty$,\\

$(ii)$ $X_n/\phi(n) \to 0$ $a.s.$ ($n \to \infty$),\\

$(iii)$ $(X_n-m_n)\rightarrow 0$ $a.s.$ $(V,1,q_n,u_n)$ for some (all) $(u_n,q_n)\in\Phi_V(\phi)$,\\

$(iv)$ $(X_n-m_n)\rightarrow 0$ $a.s.$ $(V,v_n,q_n,u_n)$ for some (all) $(u_n,v_nq_n)\in\widetilde{\Phi}_V(\phi)$.\\

\ni(b) For a given $u_n$, let $D_u(x):=\sum_{n=0}^\infty v_nx^n$ have radius of convergence $R_u\in(0,\infty)$. Let $\Phi_{uq}(\phi)$ denote the set of pairs $(u_n,q_n)$ such that $u_n\rightarrow\infty$ as $n\rightarrow\infty$, and there exists a function $h_{uq}:(0,\infty)\rightarrow(0,\infty)$ such that:
\begin{eqnarray}
\lim_{x\rightarrow\infty}h_{uq}(x)=R_u^{-1},\quad and \quad D_u(h_{uq}(n))/q_nh^n_{uq}(n)=\phi(n).\nonumber
\end{eqnarray}

\ni If $\phi^{\leftarrow}$ is subadditive, then the following two statements are equivalent:\\

$(i)$ $\mathds{E}\left[\phi^{\leftarrow}(|X|)\right] < \infty$,\\

$(ii)$ $(X_n-m_n)\rightarrow 0$ $a.s.$ $(\mathcal{P},1,q_n,v_n)$ for some (all) $(u_n,q_n)\in\Phi_V(\phi)\cap\Phi_{uq}(\phi)$.\\

\ni(c) Let $\Phi_D(\phi):=\{(u_n,q_n)\in\Phi_V(\phi):\quad u_n\in\Lambda\}$. The following two statements are equivalent:\\

$(i)$ $\mathds{E}\left[\phi^{\leftarrow}(|X|)\right] < \infty$,\\

$(ii)$ $(X_n-m_n)\rightarrow 0$ $a.s.$ $(\mathcal{V},1,q_n,u_n,\lambda)$ for some (all) $(u_n,q_n,\lambda)\in\Phi_D(\phi)\times(1,\infty)$.\\

\ni(d) If $\phi$ is regularly varying of index $\rho>0$, then the following three statements are equivalent:\\

$(i)$ $\mathds{E}\left[\phi^{\leftarrow}(|X|)\right] < \infty$,\\

$(ii)$ $\sum_1^\infty n^{-1}\mathds{P}[|\sum_{1\leq i\leq n}(X_i-m_{i+n/(\gamma-1)})|>\phi(n/(\gamma-1))\epsilon]<\infty$ $\forall\epsilon>0$ and $\forall \gamma>1$,\\

$(iii)$ $\sum_1^\infty n^{-1}\mathds{P}[\max_{1\leq k\leq n}|\sum_{1\leq i\leq k}(X_i-m_{i+n/(\gamma-1)})|>\phi(n/(\gamma-1))\epsilon]<\infty$ $\forall\epsilon>0$ and $\forall \gamma>1$.
\end{theorem}

\section{Proofs}\label{pro}

\ni{\it Proof of Theorem \ref{Conv}.} ({\it Sufficiency}) Let $\sum_{n=0}^\infty b_n/u_n$ converge. Then, by Kronecker's lemma (see, for example,~\cite{Knopp} page 129):
\begin{eqnarray}
\frac{1}{u_n}\sum_{k=0}^nb_k\rightarrow 0\quad (n\rightarrow\infty).\nonumber
\end{eqnarray}
If $a_n\rightarrow s$ as $n\rightarrow\infty$ and (\ref{conv_cond}) holds, then
\begin{eqnarray}
\frac{1}{u_n}\sum_{k=0}^n(p\circ qs)_k=\frac{1}{u_n}\sum_{k=0}^n(v_ka_k+b_k)\rightarrow s\quad(n\rightarrow\infty).\nonumber
\end{eqnarray}
\ni({\it Necessity}) Let $s_n\rightarrow s$ $(V,p_n,q_n,u_n)$. From (\ref{t2}) we have:
\begin{eqnarray}
(p\circ qs)_n=t_nu_n-t_{n-1}u_{n-1}=v_nt_{n-1}+u_n(t_n-t_{n-1}).\label{decomp}
\end{eqnarray}
If $a_n:=t_{n-1}$ and $b_n:=u_n(t_n-t_{n-1})$, then (\ref{decomp}) is the required decomposition of $(p\circ qs)_n$, since $a_n\rightarrow s$ and $\sum_{n=0}^\infty b_n/u_n$ converges.\hfill{$\Box$}\\

\ni{\it Proof of Theorem \ref{lim}.} Let $s_n\rightarrow s$ $(V,p_n,q_n,u_n)$. From (\ref{t2}) we have:
\begin{eqnarray}
(p\circ qs)_n&=&t_nu_n-t_{n-1}u_{n-1}\nonumber\\
&=&s(u_n-u_{n-1})+(t_n-s)u_n+(t_{n-1}-s)u_{n-1}=sv_n+o(u_n).\nonumber
\end{eqnarray}
\hfill{$\Box$}\\

\ni{\it Proof of Theorem \ref{oco}.} We adapt the approach of~\cite{MS}, and prove part (i) only, as the proof of part (ii) follows the same steps.\\
\ni({\it Necessity}) Let $s_n\rightarrow s$ $(V,p_n,q_n,u_n)$. For any $\alpha\in U_q$ we have:
\begin{eqnarray}
\tau_n&:=&\frac{1}{q_{\alpha(n)}-q_n}\sum_{k=n+1}^{\alpha(n)}(h\circ ut)_k=\frac{q_{\alpha(n)}s_{\alpha(n)}-q_ns_n}{q_{\alpha(n)}-q_n}\nonumber\\
&=&s_{\alpha(n)}+\frac{1}{\frac{q_{\alpha(n)}}{q_n}-1}(s_{\alpha(n)}-s_n)\rightarrow s\quad(n\rightarrow\infty).\nonumber
\end{eqnarray}
It now follows that condition (\ref{con1}) must hold:
\begin{eqnarray}
\lim_{n\rightarrow\infty}\frac{1}{q_{\alpha(n)}-q_n}\sum_{k=n+1}^{\alpha(n)}[(h\circ ut)_k-m_kt_n]=\lim_{n\rightarrow\infty}(\tau_n-t_n)=0.\nonumber
\end{eqnarray}
Similarly, for any $\beta\in L_q$ we have:
\begin{eqnarray}
\rho_n&:=&\frac{1}{q_n-q_{\beta(n)}}\sum_{k=\beta(n)+1}^{n}(h\circ ut)_k=\frac{q_ns_n-q_{\beta(n)}s_{\beta(n)}}{q_n-q_{\beta(n)}}\nonumber\\
&=&s_n+\frac{1}{\frac{q_n}{q_{\beta(n)}}-1}(s_n-s_{\beta(n)})\rightarrow s\quad(n\rightarrow\infty).\nonumber
\end{eqnarray}
It now follows that condition (\ref{con2}) must hold:
\begin{eqnarray}
\lim_{n\rightarrow\infty}\frac{1}{q_n-q_{\beta(n)}}\sum_{k=\beta(n)+1}^{n}[m_kt_n-(h\circ ut)_k]=\lim_{n\rightarrow\infty}(t_n-\rho_n)=0.\nonumber
\end{eqnarray}
\ni({\it Sufficiency}) Let the conditions (\ref{con1}) and (\ref{con2}) hold. For $\varepsilon >0$, there exists $\alpha\in U_q$ and $\beta\in L_q$ such that:
\begin{eqnarray}
-\varepsilon&\leq&\liminf_{n\rightarrow\infty}\frac{1}{q_{\alpha(n)}-q_n}\sum_{k=n+1}^{\alpha(n)}[(h\circ ut)_k-m_kt_n]\nonumber\\
\nonumber\\
&=&\liminf_{n\rightarrow\infty}(\tau_n-t_n)=s-\limsup_{n\rightarrow\infty}t_n,\nonumber\\
-\varepsilon&\leq&\liminf_{n\rightarrow\infty}\frac{1}{q_n-q_{\beta(n)}}\sum_{k=\beta(n)+1}^{n}[m_kt_n-(h\circ ut)_k]\nonumber\\
\nonumber\\
&=&\liminf_{n\rightarrow\infty}(t_n-\rho_n)=\liminf_{n\rightarrow\infty}t_n-s,\nonumber
\end{eqnarray}
which together imply $t_n\rightarrow s$ as $n\rightarrow\infty$. \hfill{$\Box$}\\

\ni{\it Proof of Theorem \ref{k}.} Let $s_n\rightarrow s$ $(V,1,q_n,u_n)$, i.e.
\begin{eqnarray}
t_n=\frac{1}{u_n}\sum_{k=0}^ns_kq_k\rightarrow s\quad (n\rightarrow\infty).\nonumber
\end{eqnarray}
We can express the sequence $s_n$ in terms of $t_n$ as:
\begin{eqnarray}
s_0=\frac{u_0}{q_0}t_0,\quad s_n=\frac{u_n}{q_n}(t_n-u_{n-1}t_{n-1}/u_{n})\quad\mbox{for $n\geq1$}.\nonumber
\end{eqnarray}
As $u_{n+1}/u_n\rightarrow 1$, we have $\hat{t}_n:=t_n-u_{n-1}t_{n-1}/u_{n}\rightarrow 0$ as $n\rightarrow\infty$. The sequence
\begin{eqnarray}
\tilde{t}_n:=\frac{1}{\tilde{u}_n}\sum_{k=0}^ns_k\tilde{q}_k=\frac{1}{\tilde{u}_n}\sum_{k=0}^n\frac{u_k}{q_k}\tilde{q}_k\hat{t}_k,\nonumber
\end{eqnarray}
is a linear transformation of the converging sequence $\hat{t}_n$. Moreover, due to assumptions  $\tilde{u}_n\rightarrow\infty$ as $n\rightarrow\infty$, and $u_n/q_n\rightarrow 1$ $(V,1,\tilde{q}_n,\tilde{u}_n)$, it is a regular transformation, and hence the conclusion.\hfill{$\Box$}\\

\ni{\it Proof of Theorem \ref{ND}.} Here we follow closely the approach of~\cite{BiG1}. To prove $(V,p_n,q_n,u_n)\Rightarrow(\mathcal{V},p_n,q_n,u_n,\lambda)$, let $t_n\rightarrow s$ $(\Omega)$. It is clear from (\ref{mov}) that
\begin{eqnarray}
c_n=t_n-\frac{u_{[w_\lambda(n)]}}{u_n}t_{[w_\lambda(n)]}.\label{map}
\end{eqnarray}
Thus, the sequence $c_n$ is a {\it transformation} of the sequence $t_n$. For each $n$, the only nonzero coefficients of such a transformation are $1$ and $u_{[w_\lambda(n)]}/u_n$. The sum of their absolute values is finite for each $n$, they shift with $n$, and their sum tends to $1-\lambda^{-1}$ as $n\rightarrow\infty$. Hence it is a regular transformation.\\

\ni To prove $(V,p_n,q_n,u_n)\Leftarrow(\mathcal{V},p_n,q_n,u_n,\lambda)$, let $c_n\rightarrow s$ $(\Omega)$. From (\ref{map}) it is clear that we can write $t_n$ as the following {\it finite} sum:
\begin{eqnarray}
t_n&=&c_n+\frac{u_{[w_\lambda(n)]}}{u_{n}}t_{[w_\lambda(n)]}\nonumber\\
&=&c_n+\frac{u_{[w_\lambda(n)]}}{u_{n}}\left[c_{[w_\lambda(n)]}+\frac{u_{[w_\lambda([w_\lambda(n)])]}}{u_{[w_\lambda (n)]}}t_{[w_\lambda([w_\lambda(n)])]}\right]\nonumber\\
&=&c_n+\frac{u_{[w_\lambda(n)]}}{u_n}c_{[w_\lambda(n)]}+\frac{u_{[w_\lambda([w_\lambda(n)])]}}{u_n}c_{[w_\lambda([w_\lambda(n)])]}+....\nonumber
\end{eqnarray}
Thus, the sequence $t_n$ can be seen as a transformation of the sequence $c_n$ with a finite number of nonzero terms. Since these coefficients are either zero or tend to zero with $n$, and their sum as $n\rightarrow\infty$ is $1+\lambda^{-1}+\lambda^{-2}+...=(1-\lambda^{-1})^{-1}$, we conclude that it is a regular transformation.\hfill{$\Box$}\\

\ni{\it Proof of Theorem \ref{Uni}.} Let (\ref{mov}) hold for all $\lambda>1$. We can write (\ref{mov}) as
\begin{eqnarray}
\frac{U(n)-U(w_\lambda(n))}{u_n}\rightarrow(1-\lambda^{-1})s\quad (n\rightarrow\infty), \quad \forall \lambda>1,\label{M1}
\end{eqnarray}
which holds also for $\lambda=1$. Define $\alpha_n:=\lambda^{-1}u_n$ and $\widetilde{U}(x) := U(u^{\leftarrow}(x))$, and rewrite (\ref{M1}) as
\begin{eqnarray}
\frac{\widetilde{U}(\lambda \alpha_n)-\widetilde{U}(\alpha_n)}{\alpha_n}\rightarrow(\lambda-1)s\quad (\alpha_n\rightarrow n),\quad \forall \lambda\geq1.\label{V}
\end{eqnarray}
Since the linear function $x$ is regularly varying of index $1$, the function $\widetilde{U}$ belongs to the de Haan class $\Pi_1$ (see Chapter 3 of~\cite{BinGT}). Hence the proof of the local uniformity follows from the proof of Theorem 3.1.16 of~\cite{BinGT} by using $\alpha_n$ instead of a continuous variable.\\

\ni{\it Proof of Theorem \ref{NDJM}}. From the previous proof it is clear that $s_n\rightarrow s$ $(\mathcal{V}_x,p_n,q_n,u(x),\lambda)$ means
\begin{eqnarray}
\frac{\widetilde{U}(\lambda y)-\widetilde{U}(y)}{y}\rightarrow(\lambda-1)s\quad (y\rightarrow\infty),\quad \forall \lambda\geq1,\label{V2}
\end{eqnarray}
where $y:=\lambda^{-1}u(x)$. From Theorem 3.2.7 of~\cite{BinGT} it now follows that (\ref{V2}) holds if and only if $s_n\rightarrow s$ $(V_x,p_n,q_n,u(x))$.\hfill{$\Box$}\\

\ni{\it Proof of Theorem \ref{UN}.} This follows from (\ref{V2}) and Theorem 3.8.4 of~\cite{BinGT}.\hfill{$\Box$}\\

\ni{\it Proof of Corollary \ref{xn}.} From Theorem \ref{NDJM},  Theorem \ref{Uni}, and Theorem \ref{ND}, respectively, it follows that:
\begin{eqnarray}
(V_x,p_n,q_n,u(x))\Leftrightarrow(\mathcal{V}_x,p_n,q_n,u(x),\lambda)\Leftrightarrow(\mathcal{V},p_n,q_n,u_n,\lambda)\Leftrightarrow(V,p_n,q_n,u_n).\nonumber
\end{eqnarray}\hfill{$\Box$}\\

\ni{\it Proof of Theorem \ref{t1}.} The following two {\it equivalence} relations are evident from the definitions of Voronoi mean and Voronoi power series:
\begin{eqnarray}
s_n\rightarrow s\quad(V,p_n,q_n,u_n)\quad\Leftrightarrow\quad (p\circ sq)_n/v_n\rightarrow s\quad (V,1,v_n,u_n),\label{equ1}\\
s_n\rightarrow s\quad(\mathcal{P},p_n,q_n,v_n)\quad\Leftrightarrow\quad (p\circ sq)_n/v_n\rightarrow s\quad (\mathcal{P},1,v_n,v_n).\label{equ2}
\end{eqnarray}
\ni(i) Here we follow closely~\cite{Ish6} and~\cite{soni}. The Voronoi power series can be written as
\begin{eqnarray}
\frac{(1-x)^{-1}\sum_{n=0}^\infty u_nt_nx^n}{(1-x)^{-1}\sum_{n=0}^\infty u_n x^n}=\frac{\sum_{n=0}^\infty u_nt_nR^n(x/R)^n}{\sum_{n=0}^\infty u_nR^n (x/R)^n}.\label{vpp}
\end{eqnarray}
If $s_n\rightarrow s$ $(V,p_n,q_n,u_n)$, then from Theorem 57 of~\cite{Har} it follows that (\ref{vpp}) converges to $s$ as $x\rightarrow R-$.\\

\ni(ii) Under the stated assumptions, it follows from Theorem 4.1 of~\cite{jakim} that $(p\circ sq)_n/v_n\rightarrow s\quad (\mathcal{P},1,v_n,v_n)$ implies $(p\circ sq)_n/v_n\rightarrow s\quad (V,1,v_n,u_n)$. The conclusion now follows from (\ref{equ2}) and (\ref{equ1}).\\

\ni(iii) The conclusions are immediate from Theorem 1.7.6 of~\cite{BinGT} and the fact that $\hat{U}(s)=(1-e^{-s})\sum_{n=0}^\infty (p\circ qs)_ne^{-ns}$.\\

\ni(iv)  Under the stated assumptions, it follows from Theorem 5.3 of~\cite{jakim} that $(p\circ sq)_n/v_n\rightarrow s\quad (\mathcal{P},1,v_n,v_n)$ implies $(p\circ sq)_n/v_n\rightarrow s$. The conclusion now follows from (\ref{equ2}).\\

\ni(v)  Under the stated assumptions, from Theorem 1 of~\cite{Ish7} we have that $(p\circ sq)_n/v_n\rightarrow s\quad (\mathcal{P},1,v_n,v_n)$ implies $(p\circ sq)_n/v_n\rightarrow s$. The conclusion now follows from (\ref{equ2}).\hfill{$\Box$}\\

\ni{\it Proof of Theorem \ref{LLN}.} $(a)$ The equivalence $(i)\Leftrightarrow (ii)$ is implicit in the proof of Theorem in~\cite{jajte}, whereas $(i)\Leftrightarrow(iii)$ follows from that theorem. The equivalence $(iii)\Leftrightarrow(iv)$ follows from Theorem 3 of~\cite{Ish}, from which we know that
\begin{eqnarray}
(V,v_n,q_n,u_n)\Leftrightarrow(V,1,v_nq_n,u_n).\nonumber
\end{eqnarray}

\ni$(b)$ Here we closely follow~\cite{BiG1}. From part $(a)$ and the Abelian result of Theorem \ref{t1} $(i)$, we have that $(i)\Rightarrow (V,1,q_n,u_n)\Rightarrow(ii)$. To prove the opposite, note that $(ii)$  implies:
\begin{eqnarray}
\frac{1}{D_u(h_{uq}(m))}\sum_{k=1}^\infty X_k^s q_k h_{uq}^k(m)=0\quad (m\rightarrow\infty)\quad a.s.,\nonumber
\end{eqnarray}
where $X_k^s=X_k-X_k'$, and $\{X_n\}_{n=1}^\infty$ and $\{X'_n\}_{n=1}^\infty$ are i.i.d.. We define
\begin{eqnarray}
\widetilde{X}_m:=\frac{1}{D_u(h_{uq}(m))}\sum_{k=1}^m X_k^s q_k h_{uq}^k(m), \quad \widehat{X}_m:=\frac{1}{D_u(h_{uq}(m))}\sum_{k=m+1}^\infty X_k^s q_k h_{uq}^k(m).\nonumber
\end{eqnarray}
Then $\widetilde{X}_m+\widehat{X}_m\rightarrow 0$ a.s., so in probability. As they are independent and symmetric, from the L\'{e}vy inequality (Lemma 2 in V.5 of~\cite{Fel}), $\widehat{X}_m\rightarrow 0$ in probability. Since $(\widetilde{X}_1,...,\widetilde{X}_m)$ and $\widehat{X}_m$ are independent, Lemma 3 of~\cite{CL} gives $\widetilde{X}_m\rightarrow 0$, $a.s.$. Repeating the same argument for
\begin{eqnarray}
\widetilde{X}_m=\frac{1}{D_u(h_{uq}(m))}\sum_{k=1}^{m-1} X_k^s q_k h_{uq}^k(m)+\frac{1}{D_u(h_{uq}(m))}X_m^s q_m h_{uq}^m(m),\nonumber
\end{eqnarray}
gives $X_m^s/\phi(m)\rightarrow 0$ $(m\rightarrow\infty)$  $a.s.$. By the Borel-Cantelli lemma, and the weak symmetrisation inequalities (pp. 257 of~\cite{Loe}),
\begin{eqnarray}
\frac{1}{2}\sum_{k=1}^\infty\mathds{P}\left[\phi^{\leftarrow}(|X-\mu_x|)\geq k\right]=\frac{1}{2}\sum_{k=1}^\infty\mathds{P}\left[|X-\mu_x|\geq \phi(k)\right]\nonumber\\
\leq\sum_{k=1}^\infty\mathds{P}\left[|X^s|\geq \phi(k)\right]<\infty,\nonumber
\end{eqnarray}
with $\mu_x$ the median of $X$, and $X^s=X-X'$, with $X$ and $X'$ i.i.d. Since $\phi^\leftarrow$ is assumed subadditive, we finally obtain:
\begin{eqnarray}
\mathds{E}\left[\phi^{\leftarrow}(|X|)\right]\leq\mathds{E}\left[\phi^{\leftarrow}(|X-\mu_x|+|\mu_x|)\right]\leq \phi^{\leftarrow}(|\mu_x|)+\mathds{E}\left[\phi^{\leftarrow}(|X-\mu_x|)\right]<\infty.\nonumber
\end{eqnarray}

\ni$(c)$ This follows immediately from part $(a)$ and Theorem \ref{ND}.\\

\ni$(d)$ Part $(a)$ and Corollary \ref{xn} show that $(i)$ is equivalent with
\begin{eqnarray}
\frac{1}{\phi(x)}\sum_{0<i\leq x}(X_i-m_i)\rightarrow 0\quad a.s.\quad (x\rightarrow\infty).\nonumber
\end{eqnarray}
By Theorem 3.2.7 of~\cite{BinGT} this is equivalent to
\begin{eqnarray}
\frac{1}{\phi(x)}\sum_{x<i\leq \gamma x}(X_i-m_i)\rightarrow 0\quad a.s.\quad (x\rightarrow\infty)\quad \forall\gamma>1.\nonumber
\end{eqnarray}
The remainder of the proof proceeds identically to that on page 1787 of~\cite{BiG1}, and is thus omitted.\hfill{$\Box$}

\section{Further remarks}\label{rem}
We give a brief account of the non-regular summability methods that appear in probability theory, analysis, and number theory.
\subsection{LLN}
As already mentioned in the introduction, the Chow-Lai laws of large numbers (LLNs) in~\cite{CL} are not regular. Further results of the same kind were also given by Li et al.~\cite{LiRJW} (double sequences of random variables). Similarly, the Marcinkiewicz-Zygmund LLN (\cite{MarZ};~\cite{Gut} \S 6.7);~\cite{Bin9} \S3) gives a non-regular summability method for $L_p$ $(0<p<2)$ when $p\neq1$ (that is except, in the Kolmogorov case), as Jajte~\cite{jajte} remarks. Generalising this, Jajte~\cite{jajte} introduces his methods, which include both regular (e.g. Ces\`{a}ro and logarithmic) and also non-regular methods.\\

Many extensions of the Kolmogorov strong LLN (SLLN) are known, in which a.s. convergence under a summability method is tied to a moment condition -- see e.g.~\cite{BinG4},~\cite{Bin10},~\cite{BiG1} -- but here the methods are regular. The main results not included here are the Marcinkiewicz-Zygmund law (above) and the Baum-Katz law (~\cite{BauK},~\cite{Gut} \S6.11,12). This has been extensively developed by  Lai~\cite{Lai2}, who introduced the idea of $r$-{\it quick convergence} (see also~\cite{BinGo}). This is essentially probabilistic, and gives, not a summability method as such, but a convergence concept giving a probabilistic analogue of a summability method -- again non-regular.\\

\subsection{Analysis}

By a theorem of Leja~\cite{Lej}, any regular N\"{o}rlund mean sums a power series at at most countably many points outside its circle of convergence. This was extended by K. Stadtm\"{u}ller to non-regular N\"{o}rlund means; her result was developed further with Grosse-Erdmann~\cite{GE-KS}.\\

Further examples of non-regular summability methods useful in analysis arise in the theory of Fourier series.  With $s_n:=\sum_{k=0}^na_k$, write
\begin{eqnarray}
\sum a_n=s\quad \mbox{or}\quad s_n\rightarrow s\quad (R,1)\quad \mbox{for}\quad \sum_i^\infty a_n\frac{\sin nh}{nh}\rightarrow s\quad (h\downarrow 0),\nonumber\\
\sum a_n=s\quad \mbox{or}\quad s_n\rightarrow s\quad (R_1)\quad \mbox{for}\quad \frac{2}{\pi}\sum_i^\infty s_n\frac{\sin nh}{nh}\rightarrow s\quad (h\downarrow 0),\nonumber\
\end{eqnarray}
Neither method is regular, and the two are not comparable. But $(R,1)$ is {\it Fourier effective} -- sums the Fourier series of any $f$ to $f$ a.e. -- which $(R_1)$ is not: there are Fourier series summable $(R_1)$ nowhere~\cite{HarR2}.\\

The $R$ here is for the Riemann, and there are Riemann methods of higher order. If one replaces $\sin nh/(nh)$ by its square, one obtains $(R,2)$, and similarly for $(R_2)$; these methods are regular~\cite{HarR1}. These methods reduce to Abel and Ces\`{a}ro methods; see~\cite{Har} App. III, \S 12.16.\\

\subsection{Number theory}

The {\it Ingham summability method} $I$ is defined by saying that
\begin{eqnarray}
s_n\rightarrow s\quad (I)\quad\mbox{if}\quad \frac{1}{x}\sum_{n\leq x}ns_n[x/n]\rightarrow s\quad (x\rightarrow\infty).\nonumber
\end{eqnarray}
This method is not regular, but can be used, together with the Wiener-Pitt (Tauberian) theorem, to prove the Prime Number Theorem (PNT), using only the non-vanishing of $\zeta$ on the $1$-line,
\begin{eqnarray}
\zeta(1+it)\neq 0\quad (t\in\mathds{R}).\label{zeta}
\end{eqnarray}
The proof of this goes back to the first proof of the prime number theorem, and has always been recognized as that property of the Riemann zeta function which is most central in the proof of this theorem (Wiener~\cite{Wie1} IV.9,~\cite{Wie2} \S17). Indeed, (\ref{zeta}) was part of Wiener's motivation in creating his Tauberian theory. The PNT is proved via Ingham's method in Hardy~\cite{Har} \S12.11; Ingham's method is developed further in Hardy~\cite{Har} App. IV. 4, Erd\H{o}s and Segal~\cite{ErdS}.

\end{document}